\documentclass[oneside,a4paper]{amsart}

\usepackage[utf8]{inputenc}
\usepackage[english]{babel}
\usepackage{amsmath}
\usepackage{amsthm}

\usepackage{amssymb}
\usepackage{mathtools}
\usepackage{blkarray}
\usepackage[unicode=true,hidelinks,bookmarksnumbered]{hyperref}
\usepackage{enumitem}
\usepackage{tikz}
\usetikzlibrary{babel}
\usetikzlibrary{cd}
\usepackage[style=numeric,maxbibnames=10]{biblatex}
\usepackage{csquotes}

\newtheorem{theorem}{Theorem}[section]
\newtheorem{corollary}[theorem]{Corollary}
\newtheorem{proposition}[theorem]{Proposition}
\newtheorem{lemma}[theorem]{Lemma}
\theoremstyle{definition}
\newtheorem{definition}[theorem]{Definition}
\newtheorem{example}[theorem]{Example}
\newtheorem{remark}[theorem]{Remark}
\newtheorem{alg}[theorem]{Algorithm}
\numberwithin{equation}{section}
\setlist{noitemsep}
\setlist[enumerate]{label=(\arabic*), font=\upshape}

\let\Im\undefined
\DeclareMathOperator\Im{Im}
\DeclareMathOperator\Ker{Ker}
\DeclareMathOperator\Hom{Hom}
\let\hom\undefined
\DeclareMathOperator\hom{hom}
\DeclareMathOperator\supp{supp}
\DeclareMathOperator\trunc{trunc}

\begin{document}
\title[Minimal flat-injective presentations over local graded rings]{On minimal flat-injective presentations over local graded rings}
\author{Fritz Grimpen}
\address{Institute for Algebra, Geometry, Topology and their Applications (ALTA), Department of Mathematics, University of Bremen, Germany}
\email{grimpen@uni-bremen.de}
\author{Anastasios Stefanou}
\address{Institute for Algebra, Geometry, Topology and their Applications (ALTA), Department of Mathematics, University of Bremen, Germany}
\email{stefanou@uni-bremen.de}
%\keywords{Flat-injective presentations, Local graded rings, Minimality, Reduction algorithm, Matlis duality}
\begin{abstract}
    Flat-injective presentations were introduced by Miller (2020) to provide combinatorial descriptions of $\mathbb Z^n$-graded modules.
    We consider them in the setting of local graded rings~$R$, with grading over an abelian group, and give a criterion for minimality of them.
    In the special case of the polynomial ring, this criterion reduces to a family of $k$-linear equations, and we are able to give an algorithmic procedure for reduction.
    Furthermore, we provide the description of a flat-injective presentation, which can be constructed from the scalar multiplication maps of a given finitely generated $R$-module.
    Thereby, we solve the construction problem for flat-injective presentations under strong finiteness assumptions.
\end{abstract}
\maketitle

\section{Introduction}
A classical problem in graded commutative algebra is the description of a given graded module in a way that is amenable for further analysis.
The classical solution is the description in terms of generators and relations, which leads to the notion of a free resolution.
Recall that a free resolution describes a module~$M$ in terms of a complex of free modules with $M$ the zeroth homology and support in the positive degrees of the complexes.
This is amenable to an analysis by methods of homological algebra for various reasons, e.g.\ by exploiting projectivity and arising long exact sequences, and this leads to a wealth theory, which is suitable for deep results.

However, free resolutions lack good compatibility with duality of modules, e.g.\ Matlis duality and character modules.
In theory, the Matlis dual of a free resolution is known to be an injective resolution of the dual, but the converse is not true.
This creates a \enquote{defect} in the lifting of Matlis duality to the category of projective modules in that it yields a functor into the category of injective modules.

One might argue that this \enquote{defect} stems from the contravariance of duality, and does not pose a problem in general.
In topological data analysis (TDA), however, this defect is incompatible with the common interpretation of a module in terms of its births (i.e.,\ generators) and its deaths (i.e.,\ cogenerators) \cite{CarlssonZomorodian2009,ZomorodianCarlsson2005}.
That is, free resolutions describe births while injective resolutions describe deaths, but a matching between births and deaths would require a relationship between a chosen free resolution and a chosen injective resolution.

This gap between free and injective resolutions is closed by the notion of flat-injective presentations, as introduced by \textcite{Miller2020a}.
Roughly speaking, a flat-injective presentation for~$M$ is a homomorphism $\varphi\colon F \to E$ with $F$ being flat, $E$ being injective, and $\Im \varphi \cong M$, see Definition~\ref{defn:flange} and \cite[Definition~5.12]{Miller2020a}.

The main feature of flat-injective presentations is that they are preserved under Matlis duality in the following sense:
If $\varphi\colon F \to E$ is a homomorphism with the property of being a flat-injective presentation, then its Matlis dual $\varphi^\vee\colon E^\vee \to F^\vee$ is one as well \cite[Proposition 5.18]{Miller2020a}.
From a structural perspective, this is not surprising, as a flat-injective presentation of~$M$ may be regarded as the composition of a flat presentation $F \to M$ with an injective presentation $M \to E$.

Furthermore, flat-injective presentations for finitely generated modules have a description in terms of matrices, see Theorem~\ref{thm:mon-mat} and \cite[Proposition 5.17]{Miller2020a}.
This makes them a good fit for computational problems, which are common in applications such as topological data analysis.

Nevertheless, the applicability of flat-injective presentations requires constructive methods.
While Miller defines minimality as a condition on the number of indecomposable direct summands, this condition is not very useful for constructive methods, since it does not give a good leverage point for constructive reduction techniques.
This leads to the problem of finding an intrinsic notion of minimality for a flat-injective presentations, which we address in this paper.

Another problem related to constructive methods is the following.
Let a $\mathbb Z^n$-graded $k[X_1, \dotsc, X_n]$-module be given by its graded components as $k$-vector spaces and its structure maps as $k$-valued matrices.
How can a (large) flat-injective presentation, which may preferably be represented by a monomial matrix, be constructed from this data?

\subsection{Related work}
\textcite{Lenzen2024} addressed both problems by describing an algorithm that transforms a minimal free resolution to a minimal flat-injective presentation.
This shifts the problem to the problem of finding a minimal free resolution for a given module, a problem that is well-investigated in computational algebra \cite{CarlssonSinghZomorodian2009,KerberRolle2021,LaScalaStillman1998,LesnickWright2022}.
However, the problem of computing a minimal flat-injective presentation can be reduced, in principle, to the problem of computing a flat cover and an injective hull.
Thus, the prior computation of a minimal flat-injective presentation may not be strictly necessary and may pose an additional burden.

Furthermore, the problem of computing an injective hull of a finitely generated $k[X_1, \dotsc, X_n]$-module~$M$ has been adressed by \textcite{HelmMiller2005}.
This provides at least a theoretical approach to the problem of constructing a flat-injective presentation for~$M$:
Given a flat cover $\hat\varphi\colon F \to M$ and an injective hull $\check\varphi\colon M \to E$, a minimal flat-injective presentation for~$M$ is given by $\varphi = \check\varphi \circ \hat\varphi$.
However, it is not clear how a monomial matrix of $\varphi$ can be constructed, because $\check\varphi$ and $\hat\varphi$ cannot be described by monomial matrices.

\subsection{Contributions}

Our main contributions in this paper are
\begin{itemize}
    \item the characterization of minimality for flat-injective presentations in Theorem~\ref{thm:minimality2},
    \item the construction of associated flat-injective presentations for arbitrary modules in Theorem~\ref{thm:assoc}, and
    \item an algorithm for reducing a free-cofree presentation of a finitely generated and finitely supported module to a minimal one, see Algorithms~\ref{alg:reduce1} and~\ref{alg:reduce2}.
\end{itemize}

Regarding the characterization of minimality, we direct our attention to the class of local graded rings~$(R, \mathfrak m, k)$.
Our definition of minimality is divided into minimality of generators and minimality of cogenerators, see Definition~\ref{defn:minimality}, and the actual characterizations are solely provided for generator-minimality, see Theorem~\ref{thm:minimality2}.
Though, by adapting Matlis duality, this dualizes to characterizations of cogenerator-minimality.

Regarding the construction of the associated flat-injective presentation, we consider the class of local graded rings~$(R, \mathfrak m, k)$ whose homogeneous subring~$R_0$ is canonically isomorphic to the residue field~$k$.
For a given $R$-module~$M$, we obtain flat and injective presentations of $M$ and the Matlis bidual $M^{\vee\vee}$, respectively, which fuse together to a flat-injective presentation, see Definition~\ref{defn:assoc} and Theorem~\ref{thm:assoc}.

Lastly, we apply both main contributions to the problem of constructing minimal flat-injective presentations of $k[X_1, \dotsc, X_n]$-modules in Section~\ref{sect:reduction}.
We prove that the notion of minimality as well as the constructed associated flat-injective presentation are amenable by monomial matrices, see Proposition~\ref{prop:assoc-mm} and Proposition~\ref{prop:reduce}.
This leads to an algorithmic procedure to construct a minimal flat-injective presentation of a finitely generated $k[X_1, \dotsc, X_n]$-module with finite support, see Algorithms~\ref{alg:reduce1} and~\ref{alg:reduce2}.
Thus, we avoid the complication of determining a minimal free resolution first, and then discarding most of its information directly afterwards.

\section{Preliminaries}

Throughout the paper, let $G$ be an additive abelian group.

\subsection{Graded rings and modules}

We begin by recalling the notions of graded rings, modules and homomorphisms.
For a general introduction to graded ring theory, we refer to the textbooks of \textcite{Hazrat2016} and \textcite[27--39]{BrunsHerzog1998}.

\begin{definition}
    Let $G$ be an additive abelian group.
    \begin{enumerate}
        \item A \emph{$G$-graded ring} is a commutative, unital ring $R$ that admits a $\mathbb Z$-module decomposition $R = \bigoplus_{a \in G} R_a$ such that $R_a R_b \subseteq R_{a + b}$ for all $a, b \in G$.
        \item A \emph{$G$-graded $R$-module} is an $R$-module~$M$ that admits a $\mathbb Z$-module decomposition $M = \bigoplus_{a \in G} M_a$ such that $R_a M_b \subseteq M_{a + b}$ for all $a, b \in G$.
        \item For $G$-graded $R$-modules $M$ and $N$, an $R$-module homomorphism~$\varphi\colon M \to N$ is said to be \emph{$G$-graded} of degree $g \in G$ if $\varphi(M_a) \subseteq N_{a+g}$ for all $a \in G$.
        The set of $G$-graded $R$-module homomorphisms from $M$ to $N$ of degree $0$ is denoted by $\Hom_R(M, N)$.
        \item An element~$m$ of a $G$-graded module $M = \bigoplus_{a \in G} M_a$ is said to be \emph{homogeneous} if $m \in M_a$ for some $a \in G$.
        \item An ideal of a $G$-graded ring $R$ is \emph{homogeneous} if it is the kernel of a graded module homomorphism $R \to M$.
        \item A $G$-graded $R$-module~$M$ is \emph{finitely generated} if it is finitely generated as an ungraded module by homogeneous elements.
        If $M$ is finitely generated and there exists a free presentation $\eta\colon F \to M$ such that $\Ker \eta$ is finitely generated, then $M$ is said to be \emph{finitely presented}.
        \item The \emph{support} of a $G$-graded module $M = \bigoplus_{a \in G} M_a$ is the set $\supp M = \{ a \in G \mid M_a \neq 0 \}$. If $\supp M$ is finite, the graded module~$M$ is said to be \emph{finitely supported}.
    \end{enumerate}
\end{definition}

\begin{example}%
    \label{ex:st-grd}
    For each $\alpha = (\alpha_1, \dotsc, \alpha_n) \in \mathbb N^n$, we define $X^\alpha = X_1^{\alpha_1} \dotsb X_n^{\alpha_n}$ in $k[X_1, \dotsc, X_n]$.
    The polynomial ring $k[X_1, \dotsc, X_n]$, for a commutative ring $k$, is $\mathbb Z^n$-graded by $k[X_1, \dotsc, X_n] = \bigoplus_{\alpha \in \mathbb N^n} kX^{\alpha}$.
    We call this the \emph{standard grading} of the polynomial ring in $n$ indeterminates.
\end{example}

\begin{definition}[Operators of graded modules]
    Let $M$ be a $G$-graded $R$-module~$M$.
    \begin{enumerate}
        \item For $g \in G$, the \emph{shift of $M$ by $g$} is the $G$-graded $R$-module~$\Sigma^g M = \bigoplus_{a \in G} M_{a - g}$.
        \item If $N$ is another $G$-graded module, we define the $G$-graded module
        \[ \hom_R(M, N) \coloneqq \bigoplus_{a \in G} \Hom_R(\Sigma^a M, N) \text. \]
        \item For $a \in G$, let $S_a$ be the subgroup of $\bigoplus_{g \in G} M_g \otimes_{\mathbb Z} N_{a - g}$ generated by elements of the form $r m \otimes n - m \otimes rn$ with $r \in R_{h}$, $m \in M_{h'}$ and $n \in N_{h''}$ for all $h, h', h'' \in G$ such that $h + h' + h'' = a$.

        The \emph{graded tensor product} of $M$ and $N$ is the $G$-graded module $M \otimes_R N = \bigoplus_{a \in G} (M \otimes_R N)_a$ with
        \[ (M \otimes_R N)_a = (\bigoplus_{g \in G} M_g \otimes_{\mathbb Z} N_{a - g})/S_a \text. \]
    \end{enumerate}
\end{definition}

\begin{remark}
    The bifunctors~$\hom_R(\mathord-, \mathord-)$ and $(\mathord- \otimes_R \mathord-)$ admit an internal tensor-hom adjunction, which is given by
    \[ \hom_R(M \otimes_R N, L) \cong \hom_R(M, \hom_R(N, L)) \text, \]
    for graded modules~$L$, $M$, and $N$.
\end{remark}

\subsection{Local graded rings and Matlis duality}

In this paper, we mostly deal with local graded rings, i.e.\ graded rings~$R$ that have a unique maximal homogeneous ideal $\mathfrak m$.
A local graded ring is denoted as a triple $(R, \mathfrak m, k)$, where $k = R/\mathfrak m$ is the \emph{residue field} of $R$.

\begin{example}
    The polynomial ring in $n$ indeterminates with coefficients in a field~$k$ and with the standard grading is a local graded ring $(k[X_1, \dotsc, X_n], \mathfrak m, k)$ with $\mathfrak m = (X_1, \dotsc, X_n)$.
\end{example}

Under sufficient finiteness assumptions, local graded rings give rise to a duality.
The underlying idea for non-graded local rings is due to \textcite{Matlis1958}, but seamlessly translates to graded local rings.

For a graded $R$-module~$M$, we denote the graded injective hull of $M$ by $E(M)$.
Existence of $E(M)$ is guaranteed because the existence argument for non-graded modules, see \cite{EckmannSchopf1953}, translates verbatim to the graded setting.

\begin{theorem}[{\textcite{Matlis1958, GotoWatanabe1978a}\relax}]%
    \label{thm:matlis}
    Let $R$ be a Noetherian local graded ring, i.e.\ every homogeneous ideal of $R$ is finitely generated.
    \begin{enumerate}
        \item\label{thm:matlis:1} A graded $R$-module is indecomposable and graded injective if and only if it is isomorphic to $\Sigma^g E(R/\mathfrak p)$ for some homogeneous prime ideal~$\mathfrak p$ and some~$g \in G$.
        \item\label{thm:matlis:2} Every injective, graded $R$-module~$M$ admits a decomposition
        \[ M \cong \bigoplus_{i \in I} \Sigma^{g_i} E(R/\mathfrak p_i) \text, \]
        where $(\mathfrak p_i)_{i \in I}$~is a family of homogeneous prime ideals of $R$, and where $(g_i)_{i \in I}$~is a (non-unique) family in~$G$.
        The family $(\mathfrak p_i)_{i \in I}$ is uniquely determined by the module~$M$.
        \item\label{thm:matlis:4} The contravariant, exact functor $\hom_R(\mathord-, E(k))$ establishes a one-to-one correspondence from graded flat modules to graded injective modules.
    \end{enumerate}

    \begin{proof}
        Parts~\ref{thm:matlis:1} and~\ref{thm:matlis:2} can be proven by translating the arguments in \cite[Theorem 2.5 and Proposition 3.1]{Matlis1958} to graded rings, cf.\ \cite[Theorem~1.2.1]{GotoWatanabe1978a}.
        We only prove \ref{thm:matlis:4} for the reader's convenience; the same argument is given in \cite[\pno~219, Lemma~11.23]{MillerSturmfels2005}.

        If $N$ is a graded flat module, the functor $\mathord- \otimes_R N$ is exact.
        Thus, $\Hom_R(\mathord- \otimes_R N, E(k))$ is an exact functor, which is naturally isomorphic to the exact functor $\Hom_R(\mathord-, \hom_R(N, E(k)))$.
        That is, $\hom_R(N, E(k))$ is graded injective.

        If $\hom_R(N, E(k))$ is a graded injective module, the functor $\Hom_R(\mathord- \otimes_R N, E(k))$ is exact.
        Recall that $\Hom_R(\mathord-, E(k))$ reflects exactness \cite[\pno~64, Proposition 2.42]{Rotman2009}.
        Thus, the functor $\mathord- \otimes_R N$ must be exact.
        Hence, $N$ is graded flat.
    \end{proof}
\end{theorem}

\begin{definition}
    Let $(R, \mathfrak m, k)$ be a local graded ring.
    The \emph{Matlis dual} of a graded $R$-module~$M$ is defined by $M^\vee = \hom_R(M, E(k))$.
\end{definition}

\begin{remark}
    Suppose that $M$ is pointwise finite-dimensional, i.e.\ the $R_0$-module $M_a$ is finite-dimensional for every $a \in G$.
    Then, the evaluation homomorphism
    \[ u_M\colon M \to M^{\vee\vee}\text, \quad u_M(m)(f') = f'(m)\text, \]
    is an isomorphism.
\end{remark}

\begin{example}%
    \label{ex:action}
    Again, we consider the polynomial ring in $n$ indeterminates with the standard grading, see Example~\ref{ex:st-grd}.
    The (graded) injective hull~$E(k)$ of $k$ has a copy of $k$ in each degree in the negative orthant $-\mathbb N^n$.
    
    The Matlis dual of a graded $k[X_1, \dotsc, X_n]$-module~$M$ is therefore given by
    \[ M^\vee = \hom_R(M, E(k)) = \bigoplus_{g \in G} \Hom_R(\Sigma^g M, E(k)) \cong \bigoplus_{g \in G} \Hom_k(M_{-g}, k) \text. \]

    In particular, $R^\vee$ is the graded $k[X_1, \dotsc, X_n]$-module with a copy of~$k$ in each degree in the negative orthant $-\mathbb N^n$.
    The action of $R$ on $R^\vee$ given by
    \[ X^{\beta} \cdot r' = \begin{cases*}
        r' & if $a + \beta \in -\mathbb N^n$, \\
        0 & otherwise,
    \end{cases*} \]
    for $r' \in R^\vee_a$, $a \in \mathbb Z^n$, and $\beta \in \mathbb N^n$.
\end{example}

\begin{definition}[Cofree modules]%
    \label{defn:cogenerator}
    A graded module~$M$ is said to be \emph{cofree} if it is the Matlis dual of a (graded) free module.
        That is, if there exists a (graded) free module~$F$ such that $M^{\vee} = F$.
\end{definition}

\begin{remark}
    If a free module~$F$ is reflexive, i.e. the evaluation homomorphism $M \to M^{\vee\vee}$ is an isomorphism, then the associated cofree module $F^\vee$ is reflexive as well.
\end{remark}

\section{Flat-injective presentations over local graded rings}%
\label{sect:fringe}

Throughout this section, let $G$ be an additive, abelian group, and let $(R, \mathfrak m, k)$ be a local $G$-graded ring.
All modules are understood to be $G$-graded $R$-modules.

\subsection{Definitions}

Flat-injective presentations of $\mathbb Z^n$-graded modules were introduced by \textcite[Definition~5.12]{Miller2020a}.
We generalize this notion to the setting of $G$-graded $R$-modules.

\begin{definition}[Flat-injective presentation]%
    \label{defn:flange}
    Let $M$ be a $G$-graded $R$-module.
    A \emph{flat-injective presentation for $M$} is a graded homomorphism $\varphi\colon F \to E$ of degree 0 such that $F$ is graded flat, $E$ is graded injective, and $\Im \varphi \cong M$.
    
    If $F$ is a free module, then $\varphi$ is called a~\emph{free-injective presentation}, and if $E$ is additionally cofree (see Definition~\ref{defn:cogenerator}), the free-injective presentation~$\varphi$ is called a~\emph{free-cofree presentation}.
\end{definition}

\subsection{Matlis duality of flat-injective presentations}

We give a brief overview over Matlis duality of flat-injective presentations.
The general slogan is that the Matlis dual of a flat-injective presentations for~$M$ is a flat-injective presentation for~$M^\vee$.

\begin{proposition}%
    \label{prop:matlis-fi}
    Let $\varphi\colon F \to E$ be a flat-injective presentation for a module~$M$.
    \begin{enumerate}
        \item\label{prop:matlis-fi:1} The Matlis dual~$\varphi^\vee\colon E^\vee \to F^\vee$ is a flat-injective presentation for~$M^\vee$.
        \item\label{prop:matlis-fi:2} If $\varphi$ is a free-cofree presentation, then the Matlis dual~$\varphi^\vee$ is a free-cofree presentation for~$M^\vee$.
    \end{enumerate}

    \begin{proof}
        \ref{prop:matlis-fi:1}:
        Let $\hat\varphi\colon F \to M$ be an epimorphism and let $\check\varphi\colon M \to E$ be a monomorphism such that $\varphi = \check\varphi \circ \hat\varphi$.
        Because Matlis duality is exact, we have a monomorphism $(\hat\varphi)^\vee\colon M^\vee \to F^\vee$ and an epimorphism $(\check\varphi)^\vee\colon E^\vee \to M^\vee$.
        It follows that $\varphi^\vee = (\hat\varphi)^\vee \circ (\check\varphi)^\vee$.
        Hence, $\Im \varphi^\vee \cong \Im (\check\varphi)^\vee \cong M^\vee$.

        Furthermore, $F^\vee$ is graded injective and $E^\vee$ is graded flat by Theorem~\ref{thm:matlis}~\ref{thm:matlis:4}.

        \ref{prop:matlis-fi:2}:
        This follows immediately by using the definition of a free-cofree presentation and part~\ref{prop:matlis-fi:1}.
    \end{proof}
\end{proposition}

\subsection{Minimality of free-injective presentations}

\begin{definition}%
    \label{defn:minimality}
    Suppose that $\varphi\colon F \to E$ is a flat-injective presentation for a $G$-graded module $M$.
    The flat-injective presentation~$\varphi$ is said to be
    \begin{enumerate}
        \item \emph{generator-minimal} if $\Ker \varphi \subseteq \mathfrak mF$,
        \item \emph{cogenerator-minimal} if $\varphi^\vee$ is generator-minimal, and
        \item \emph{minimal} if it is generator-minimal and cogenerator-minimal.
    \end{enumerate}
\end{definition}

We immediately obtain the following characterization of generator-minimal free-injective presentations.

\begin{theorem}%
    \label{thm:minimality2}
    Let $M$ be a module, and suppose that
    \[ \varphi\colon \bigoplus_{j \in J} \Sigma^{\beta_j} R \longrightarrow E \]
    is a free-injective presentation for~$M$.
    For $j \in J$, denote by $e_j$ the $j$-th canonical basis element of $\bigoplus_{j \in J} \Sigma^{\beta_j} R$ in degree~$\beta_j \in G$.
    The following assertions are equivalent.
    \begin{enumerate}[label=(\roman*)]
        \item\label{thm:minimality2:1}
        The free-injective presentation~$\varphi$ is generator-minimal.
        \item\label{thm:minimality2:2}
        The graded epimorphism~$\varphi\colon \bigoplus_{j \in J} \Sigma^{\beta_j} R \to M$ is a projective cover of~$M$.
        \item\label{thm:minimality2:3}
        The epimorphism~$\varphi\colon \bigoplus_{j \in J} \Sigma^{\beta_j} R \to M$ extends to a minimal free resolution of $G$-graded $R$-modules of~$M$, i.e.\ there is an exact sequence
        \[ \begin{tikzcd}
            \cdots \rar & F_2 \rar["\partial_2"] & F_1 \rar["\partial_1"] & \bigoplus_{j \in J} \Sigma^{\beta_j} R \rar["\varphi"] & M \rar & 0
        \end{tikzcd} \]
        with $F_1, F_2, \dotsc$ being graded free $R$-modules such that
        \[ \Im \partial_i \subseteq \mathfrak m F_{i-1} \qquad \text{for $i = 1, 2, \dotsc$,} \]
        where $F_0 = \bigoplus_{j \in J} \Sigma^{\beta_j} R$.
        \item\label{thm:minimality2:4}
        It holds $\Ker \varphi \subseteq \bigoplus_{j \in J} \mathfrak m \Sigma^{\beta_j} R$.
        \item\label{thm:minimality2:5}
        For all $\sigma_1, \dotsc, \sigma_s \in R$ and distinct $j_1, \dotsc, j_s \in J$, $\varphi(\sigma_1 e_{j_1} + \dotsb + \sigma_s e_{j_s}) = 0$ implies
        \[ \sigma_1 e_{j_1} + \dotsb + \sigma_s e_{j_s} \in \bigoplus_{j \in J} \mathfrak m \Sigma^{\beta_j} R \text. \]
        \item\label{thm:minimality2:6}
        For all $\sigma_1, \dotsc, \sigma_s \in R$ and distinct $j_1, \dotsc, j_s \in J$, $\varphi(\sigma_1 e_{j_1} + \dotsb + \sigma_s e_{j_s}) = 0$ implies $\sigma_1, \dotsc, \sigma_s \in \mathfrak m$.
    \end{enumerate}

    \begin{proof}
        \ref{thm:minimality2:1}, \ref{thm:minimality2:2}, \ref{thm:minimality2:3}, and \ref{thm:minimality2:4} are equivalent by standard arguments, cf.\ \cite{Eisenbud1995}.
        \ref{thm:minimality2:4} to \ref{thm:minimality2:5} is trivial.

        We consider \ref{thm:minimality2:5} to \ref{thm:minimality2:4}:
        Let $m \in \Ker \varphi \subseteq \bigoplus_{j \in J} \Sigma^{\beta_j} R$.
        That is, there exist $\sigma_1, \dotsc, \sigma_s \in R$ and distinct $j_1, \dotsc, j_s \in J$ such that $\varphi(\sigma_1 e_{j_1} + \dotsb + \sigma_s e_{j_s}) = 0$.
        Thus, we have
        \[ m = \sigma_1 e_{j_1} + \dotsb + \sigma_s e_{j_s} \in \bigoplus_{j \in J} \mathfrak m\Sigma^{\beta_j} R \text, \]
        where the second equality holds by assumption.

        \ref{thm:minimality2:5} to \ref{thm:minimality2:6}:
        Let $\sigma_1, \dotsc, \sigma_s \in R$ and $j_1, \dotsc, j_s \in J$ distinct such that $m \coloneqq \sigma_1 b'_{j_1} + \dotsb + \sigma_s b'_{j_s} \in \Ker \varphi$.
        By assumption we have $m \in \bigoplus_{j \in J} \mathfrak m \Sigma^{\beta_j} R$.
        Moreover, we may assume without loss of generality that $m$ is homogeneous since $\varphi$ is graded.

        Then, since the elements $e_{j_1}, \dotsc, e_{j_s}$ are $R$-linearly independent, it follows that $\sigma_1, \dotsc, \sigma_s \in \mathfrak m$.
    \end{proof}
\end{theorem}

\begin{remark}
    The characterization of generator-minimality in Theorem~\ref{thm:minimality2}~\ref{thm:minimality2:3} relates it to typical minimality of a free resolution, cf.\ \cite[472--473]{Eisenbud1995}.
\end{remark}

As a corollary, we obtain that, in the assertions Theorem~\ref{thm:minimality2}~\ref{thm:minimality2:5} and~\ref{thm:minimality2:6}, it suffices to assume that $\sigma_1, \dotsc, \sigma_s \in R$ are homogeneous elements.
The argument essentially uses that $\varphi\colon F \to E$ is a graded homomorphism.

\begin{corollary}
    \label{cor:minimality}
    With the notation of Theorem~\ref{thm:minimality2}, the following assertions are equivalent:
    \begin{enumerate}[label=(\roman*)]
        \item\label{cor:minimality:1} The free-injective presentation $\varphi$ is generator-minimal.
        \item\label{cor:minimality:2} For all homogeneous $\sigma_1, \dotsc, \sigma_s \in R$ and distinct $j_1, \dotsc, j_s \in J$, $\varphi(\sigma_1 e_{j_1} + \dotsb + \sigma_s e_{j_s}) = 0$ implies
        \[ \sigma_1 e_{j_1} + \dotsb + \sigma_s e_{j_s} \in \bigoplus_{j \in J} \mathfrak m \Sigma^{\beta_j} R \text. \]
        \item\label{cor:minimality:3} For all homogeneous $\sigma_1, \dotsc, \sigma_s \in R$ and distinct $j_1, \dotsc, j_s \in J$, $\varphi(\sigma_1e_{j_1} + \dotsb + \sigma_se_{j_s}) = 0$ implies $\sigma_1, \dotsc, \sigma_s \in \mathfrak m$.
    \end{enumerate}

    \begin{proof}
        The implications \ref{cor:minimality:1} to \ref{cor:minimality:2} to \ref{cor:minimality:3} follow directly from Theorem~\ref{thm:minimality2}.
        We prove that the homogeneity assumption on $\sigma_1, \dotsc, \sigma_s \in R$ in \ref{cor:minimality:3} is superfluous.
        Then, \ref{cor:minimality:3} to \ref{cor:minimality:1} follows.

        Let $\sigma_1, \dotsc, \sigma_s \in R$ and let $j_1, \dotsc, j_s \in J$ be distinct.
        Suppose that \ref{cor:minimality:3} holds, and assume that
        \[ \varphi(\sigma_1 e_{j_1} + \dotsb + \sigma_s e_{j_s}) = 0 \text. \]

        For each $a \in G$, we have
        \[ (\sigma_1 \varphi(e_{j_1}))_a + \dotsb + (\sigma_s \varphi(e_{j_s}))_a = 0 \text. \]
        Since the homomorphism $\varphi\colon F \to E$ is graded of degree 0, the elements $\varphi(e_{j_l})$, where $l \in \{ 1, \dotsc, s \}$, are homogeneous in $E$.
        
        For each $l \in \{ 1, \dotsc, s \}$, we can write $\sigma_s \varphi(e_{j_l}) = \sum_{g \in G} (\sigma_l)_g \varphi(e_{j_l})$, where for $g \in G$ the summand $(\sigma_l)_g \varphi(e_{j_l})$ is contained in $E_{\beta_{j_l} + a}$.
        In other words, the scalar multiplication map $(\mathord-\cdot\varphi(e_{j_l}))\colon R \to E$ is graded of degree $\beta_{j_l}$.

        Since $G$ is a group, the unique~$g \in G$ such that $(\sigma_s)_g \varphi(e_{j_l}) \in E_a$ is given by $g = a - \beta_{j_l}$.
        Therefore,
        \[ (\sigma_l \varphi(e_{j_l}))_a = (\sigma_l)_{\beta_{j_l} - a} \varphi(e_{j_l}) \text, \]
        for all $l \in \{ 1, \dotsc, s \}$.

        Now, applying \ref{cor:minimality:3} to
        \[ (\sigma_1)_{\beta_{j_1} - a} \varphi(e_{j_1}) + \dotsb + (\sigma_s)_{\beta_{j_s} - a} \varphi(e_{j_s}) = 0 \]
        for every $a \in G$ gives that $\sigma_1, \dotsc, \sigma_s \in \mathfrak m$.
    \end{proof}
\end{corollary}

\subsection{Associated flat-injective presentation}%
\label{sect:fringe:assoc}

We assume that every local graded ring~$(R, \mathfrak m, k)$ satisfies the property that the canonical monomorphism~$R_0 \to k$, $r \mapsto r + \mathfrak m$, is an isomorphism.
That is, the subring~$R_0$ of~$R$ is a field.

A standard example of a local graded ring, which satisfies this assumption, is the polynomial ring~$R = k[X_1, \dotsc, X_n]$ in~$n$ indeterminates with the standard grading.

A counterexample is the ring of complex numbers~$\mathbb C = \mathbb R \oplus \mathbb R$ with the obvious $\mathbb Z_2$-grading $\mathbb C = \mathbb R \oplus i\mathbb R$.

\begin{definition}[Associated flat-injective presentation]%
    \label{defn:assoc}
    Let $M = \bigoplus_{a \in G} M_a$ be a graded module.
    \begin{enumerate}
        \item The \emph{associated flat module}~$F_M$ of $M$ is defined by
        \[ F_M = \bigoplus_{g \in G} \Sigma^g R \otimes_k M_g \text. \]
        \item The \emph{associated injective module}~$E_M$ of $M$ is defined by
        \[ E_M = \prod_{h \in G} \Hom_k(M^\vee_h, \Sigma^{-h} R^\vee) \text. \]
        \item Let $\hat\varphi_M\colon F_M \to M$ be the graded homomorphism induced by
        \begin{align*}
            \Sigma^g R \otimes_k M_g &\longrightarrow M \text, \\
            r \otimes m &\longmapsto rm \text.
        \end{align*}
        Further, let $\check\varphi_M\colon M^{\vee\vee} \to E_M$ be the graded homomorphism induced by
        \begin{align*}
            M^{\vee\vee} &\longrightarrow \Hom_k(M^\vee_h, \Sigma^{-h} R^\vee) \text, \\
            f'' &\longmapsto (f' \mapsto (r \mapsto r \cdot f''(f'))) \text,
        \end{align*}
        for $r \in \Sigma^{-h} R$, $f'' \in M^{\vee\vee}$, and $f' \in M^\vee_h$.
        The \emph{associated flat-injective presentation} of $M$ is the homomorphism
        \[ \varphi_M = \check\varphi_M \circ u_M \circ \hat\varphi_M\colon F_M \to E_M \text, \]
        where $u_M\colon M \to M^{\vee\vee}$ is the evaluation homomorphism.
    \end{enumerate}
\end{definition}

\begin{theorem}%
    \label{thm:assoc}
    The associated flat-injective presentation~$\varphi_M$ of $M$ is a flat-injective presentation~for $M$ in the sense of Definition~\ref{defn:flange}.

    \begin{proof}
        Surjectivity and injectivity of $\hat\varphi_M\colon F_M \to M$ and $\check\varphi_M\colon M^{\vee\vee} \to E_M$, respectively, follow immediately.
        The evaluation homomorphism~$u_M\colon M \to M^{\vee\vee}$ is injective because $(\mathord-)^\vee$ is a faithful functor.
        Thus, $\Im \varphi_M \cong \Im \hat\varphi_M \cong M$.

        For fixed $g \in G$ consider the graded $R$-module $\Sigma^gR \otimes_k M_g$.
        The tensor-hom adjunction gives the natural isomorphisms
        \begin{align*}
            \Hom_R(\Sigma^gR \otimes_k M_g, \mathord-) &\cong \Hom_k(M_g, \Hom_R(\Sigma^gR, \mathord-)) \text.
        \end{align*}
        The $k$-vector space $M_g$ and the graded $R$-module $\Sigma^g R$ are projective.
        Therefore, $\Sigma^g R \otimes_k M_g$ is a graded projective $R$-module.

        By a similar argument, $\Hom_k(M^\vee_h, \Sigma^{-h} R^\vee)$ is a graded injective $R$-module for every $h \in G$.

        Hence, $F_M$ and $E_M$ are graded flat and graded injective, respectively.
    \end{proof}
\end{theorem}

\begin{remark}
    Let $M$ be a graded module.
    \begin{enumerate}
        \item The module~$F_M$ is graded free since it is graded projective and $R$ is graded local, see \cite[\pno~36, Proposition~1.5.15 (d)]{BrunsHerzog1998}.
        \item The module~$E_M$ is graded cofree since $E_M \cong (F_{M^\vee})^\vee$.
    \end{enumerate}
    Thus, the associated flat-injective presentation~$\varphi_M$ of~$M$ is a free-cofree presentation.
\end{remark}

\begin{lemma}%
    \label{lem:assoc-flange}
    For $r \otimes m \in \Sigma^g R \otimes_k M_g$ it holds
    \[ \varphi_M(r \otimes m)(f') = f'(rm) \]
    for all $f' \in M^\vee$.

    \begin{proof}
        Let $u_M\colon M \to M^{\vee\vee}$ be the evaluation homomorphism, which is given by $u_M(m)(f') = f'(m)$ for $f' \in M^\vee$ and $m \in M$.

        Without loss of generality, we assume that $f'$ is homogeneous of degree zero, i.e.\ $f' \in \Hom_R(M, E(k))$.
        By substituting the definitions of $\hat\varphi_M$, $\check\varphi_M$, and $u_M$ in this order, we have
        \begin{align*}
            \varphi_M(r \otimes m)(f') &= (\check\varphi_M \circ u_M \circ \hat\varphi_M)(r \otimes m)(f') \\
            &= (\check\varphi_M \circ u_M)(rm)(f') = \check\varphi_M(u_M(rm))(f') \\
            &= u_M(rm)(f') \\
            &= f'(rm) \text.
        \end{align*}
    \end{proof}
\end{lemma}

\section{A reduction algorithm for free-injective presentations}%
\label{sect:reduction}

In this section, we consider monomial matrices as a computational tool to handle and minimize free-injective presentations of $\mathbb Z^n$-modules.
We regard the polynomial ring $R = k[X_1, \dotsc, X_n]$ with the standard grading over $G = \mathbb Z^n$.

Further, we endow $\mathbb Z^n$ with the order
\[ (a_1, \dotsc, a_n) \leq (b_1, \dotsc, b_n) \iff \text{$a_i \leq b_i$ for all $i \in \{ 1, \dotsc, n \}$,} \]
for $(a_1, \dotsc, a_n), (b_1, \dotsc, b_n) \in \mathbb Z^n$.
In the examples, we shall deliberately identify $k[X_1, \dotsc, X_n]$-modules with diagrams indexed by $\mathbb Z^n$; these diagrams are also known as \emph{persistence modules}, see\ \cite{ZomorodianCarlsson2005,CarlssonZomorodian2009}.

We denote the Laurent polynomial ring, which is the localization of $k[X_1, \dotsc, X_n]$ along $\{ X^\alpha \mid \alpha \in \mathbb N^n \}$, by $k[X_1^{\pm1}, \dotsc, X_n^{\pm1}]$.

\subsection{Monomial matrices}

A crucial observation by Miller is that free-cofree presentations and monomial matrices whose rows and columns are labelled by degrees are in one-to-one correspondence \cite{Miller2020a}.
This correspondence exposes monomial matrices as the adequate tool for manipulation of free-cofree presentations.

\begin{definition}[Monomial matrix]
    Let $(\alpha_i)_{i = 1, \dotsc, r}$ and $(\beta_j)_{j = 1, \dotsc, s}$ be families in $\mathbb Z^n$.
    \begin{enumerate}
        \item A \emph{monomial matrix} of size $(r, s)$ is a $k$-valued $(r \times s)$-matrix $A = (a_{ij})_{i, j}$ such that $a_{ij} = 0$ unless $\alpha_i \geq \beta_j$, or equivalently, unless $X^{\beta_j}$ divides $X^{\alpha_i}$ in $k[X_1^{\pm1}, \dotsc, X_n^{\pm1}]$.
        \item In the context of a monomial matrix, the elements of the family $(\alpha_i)_i$ are called \emph{cogenerator degrees}, and the elements of the family $(\beta_j)_j$ are called \emph{generator degrees}.
        \item The \emph{addition} and \emph{scalar multiplication} (over $k$) on monomial matrices, which have fixed size and (co)generator degrees, are defined entrywise.
    \end{enumerate}
\end{definition}

We write a concrete monomial matrix $A = [a_{ij}]_{i, j}$ as
\[ \begin{blockarray}{cccccc}
    & \beta_1 & \cdots & \beta_j & \cdots & \beta_s \\
    \begin{block}{c[ccccc]}
        \alpha_1\mathstrut & a_{11}\mathstrut & \cdots & a_{1j} & \cdots & a_{1s} \\
        \vdots & \vdots & & \vdots & & \vdots \\
        \alpha_i & a_{i1} & \cdots & a_{ij} & \cdots & a_{is} \\
        \vdots & \vdots & & \vdots & & \vdots \\
        \alpha_r\mathstrut & a_{r1}\mathstrut & \cdots & a_{rj} & \cdots & a_{rs} \\
    \end{block}
\end{blockarray} \text. \]

\begin{theorem}[{\textcite[Proposition 5.17]{Miller2020a}}]%
    \label{thm:mon-mat}
    Let $F = \bigoplus_{j = 1}^s \Sigma^{\beta_j} R$ and $E = \bigoplus_{i = 1}^r \Sigma^{\alpha_i} R^\vee$, where $\alpha_i, \beta_j \in G$ for $i \in \{ 1, \dotsc, r \}$ and $j \in \{ 1, \dotsc, s \}$, respectively.

    There exists a $k$-isomorphism between the graded homomorphisms $\varphi\colon F \to E$ of degree 0 and the $(r \times s)$-monomial matrices with cogenerator degrees $(\alpha_i)_{i = 1, \dotsc, r}$ and generator degrees $(\beta_j)_{j = 1, \dotsc, s}$.
\end{theorem}

We describe the monomial matrix for the associated flat-injective presentation~$\varphi_M$ of a finitely generated and finitely supported module~$M$.

\begin{proposition}%
    \label{prop:assoc-mm}
    Let $M$ be a graded module that is finitely generated and finitely supported.
    Write $d_g = \dim_k M_g$ for $g \in G$.
    Let
    \[ I = J = \bigcup_{g \in \supp M} \{ g \} \times \{ 1, \dotsc, d_g \} \text. \]
    Further, choose for every $g \in \supp M$ a basis~$(e_{(g, 1)}, \dotsc, e_{(g, d_g)})$ of the $k$-vector space $M_g$.

    Consider the monomial matrix $A = [a_{(g, \mu)(h, \nu)}]_{(g, \mu) \in I, (h, \nu) \in J}$ with cogenerator degrees $(g)_{(g, \mu) \in I}$ and generator degrees $(h)_{(h, \nu) \in J}$ that is defined as follows.
    For $(g, \mu) \in I$ and $(h, \nu) \in J$ such that $h \leq g$, let $a_{(g, \mu)(h, \nu)}$ be the $\mu$-th coefficient of $X^{g-h}\cdot e_{(h, \nu)}$ in $M_g$.
    For all other $(g, \mu) \in I$ and $(h, \nu) \in J$, let $a_{(g, \mu)(h, \nu)} = 0$.

    That is, the underlying matrix~$A$ is chosen such that
    \[ X^{g-h} \cdot e_{(h, \nu)} = a_{(g, 1)(h, \nu)} e_{(g, 1)} + \dotsb + a_{(g, d_g)(h, \nu)} e_{(g, d_g)} \]
    holds for all $(g, \mu) \in I$ and $(h, \nu) \in J$.

    Then, this monomial matrix~$A$ represents the associated free-cofree presentation of~$M$ in the sense of Theorem~\ref{thm:mon-mat}.

    \begin{proof}
        For $(g, \mu) \in I$ denote by $e^\ast_{(g, \mu)} \in \Hom_k(M_g, k)$ the dual vector that is defined by
        \[ e^\ast_{(g, \mu)}(\lambda_1 e_{(g, 1)} + \dotsb + \lambda_{d_g} e_{(g, d_g)}) = \lambda_\mu \text. \]
        Since $M$ is finitely generated, the family $(e^\ast_{(g, \mu)})_{\mu = 1, \dotsc, d_g}$ is a basis of $\Hom_k(M_g, k) \cong M_g$.

        First, we note that
        \[ F_M \cong \bigoplus_{g \in G} \Sigma^g R^{d_g} \quad \text{and} \quad E_M \cong \bigoplus_{h \in G} \Sigma^h (R^\vee)^{d_h} \text. \]
        This follows from $M$ being finitely generated and finitely supported, respectively.
        Thus, the described matrix~$A$ has the same generator and cogenerator degrees as $\varphi_M\colon F_M \to E_M$.

        We use Lemma~\ref{lem:assoc-flange} to show that the described coefficients $a_{(g, \mu)(h, \nu)}$ indeed correspond to~$\varphi_M$:
        For $(g, \mu) \in I$ and $(h, \nu) \in J$ such that $h \leq g$, we have
        \[ \varphi_M(X^{g-h} \otimes e_{(h, \nu)})(e^\ast_{(g, \mu)}) = e^\ast_{(g, \mu)} (X^{g - h} e_{(h, \nu)}) \text. \]
        Then, since $X^{g - h} e_{(h, \nu)} \in M_h$, there exist uniquely determined coefficients
        \[ a_{(g, 1)(h, \nu)}, \dotsc, a_{(g, d_g)(h, \nu)} \in k \]
        such that
        \[ X^{g-h} e_{(h, \nu)} = a_{(g, 1)(h, \nu)} e_{(g, 1)} + \dotsb + a_{(g, d_g)(h, \nu)} e_{(g, d_g)} \text. \]
        Thus, $\varphi_M(X^{g-h} \otimes e_{(h, \nu)})(e^\ast_{(g, \mu)}) = a_{(g, \mu)(h, \nu)}$.
        It follows, by linearity of $\varphi_M$ and $\varphi_M(X^{g-h} \otimes e_{(h, \nu)})$ for all $(g, \mu) \in I$ and $(h, \nu) \in J$, that this determines the matrix coefficients of~$A$.
    \end{proof}
\end{proposition}

\begin{example}
    We consider the $k[X_1]$-module~$M$ that is given by the diagram
    \[ \begin{tikzcd}[ampersand replacement=\&, row sep=0pt]
        \cdots \rar \& 0 \rar \& k^2 \rar["{\left(\begin{smallmatrix}1 & 1\end{smallmatrix}\right)}"] \& k \rar \& 0 \rar \& \cdots \\
        \& \& \scriptstyle(0) \& \scriptstyle(1) \& \&
    \end{tikzcd} \]
    and that is supported in degrees $0$ and $1$.

    The associated flat-injective presentation of~$M$ is given by
    \begin{align*}
        \varphi_M\colon \Sigma^0 R^2 \oplus \Sigma^1 R &\longrightarrow \Sigma^0 (R^\vee)^2 \oplus \Sigma^1 R^\vee \text, \\
        (X^\alpha, X^\beta, X^\gamma) &\longmapsto (\delta_{\alpha \leq 0} X^\alpha, \delta_{\beta \leq 0} X^\beta, \delta_{\gamma \leq 1} X^\gamma) \text.
    \end{align*}
    The monomial matrix of~$\varphi_M$ is given by
    \[ \begin{blockarray}{cccc}
        & X_1^0 & X_1^0 & X_1^1 \\
        \begin{block}{c[ccc]}
            X_1^0 & 1 & 0 & 0 \\
            X_1^0 & 0 & 1 & 0 \\
            X_1^1 & 1 & 1 & 1 \\
        \end{block}
    \end{blockarray} \text. \]
\end{example}

\subsection{Reduction of generators}

Purely for notational convenience, we introduce the following convention on bases of free modules.
For a given free module $F = \bigoplus_{j = 1}^s \Sigma^{\beta_j} R$, let the \emph{canonical homogeneous basis} be denoted by $\mathcal B'_F = (e_1, \dotsc, e_s)$ with $e_j = 1 \in \Sigma^{\beta_j} R$ for every $j \in \{ 1, \dotsc, s \}$.
Under this convention, the element $e_j$ is a generator of $\Sigma^{\beta_j} R$ of degree $\beta_j$.

\begin{proposition}%
    \label{prop:reduce}
    Let $(\beta_j)_{j = 1, \dotsc, s}$ be a family in $\mathbb Z^n$, and suppose that
    \[ \varphi\colon \bigoplus_{j=1}^s \Sigma^{\beta_j} R \longrightarrow E \]
    is a free-injective presentation for a module~$M$.
    \begin{enumerate}
        \item\label{prop:reduce:1} If there exist $1 \leq j_0 \leq s$ and $\lambda_{1}, \dotsc, \lambda_{j_0-1}, \lambda_{j_0+1}, \dotsc, \lambda_s \in k$ such that
        \begin{equation}%
            \label{eqn:minimality3}
            \varphi(e_{j_0}) = \sum_{\substack{j \neq j_0 \text, \\ \beta_{j_0} \geq \beta_j}} \lambda_{j} X^{\beta_{j_0} - \beta_j} \varphi(e_j) \text,
        \end{equation}
        then the restriction
        \[ \tilde\varphi_{j_0}\colon \bigoplus_{j \neq j_0} \Sigma^{\beta_j} R \longrightarrow E \]
        is a free-injective presentation for~$M$.
        \item\label{prop:reduce:2} The free-injective presentation~$\varphi$ is generator-minimal if and only if, for all $1 \leq j_0 \leq s$ and $\lambda_1, \dotsc, \lambda_{j_0-1}, \lambda_{j_0+1}, \dotsc, \lambda_{s} \in k$, it holds
        \[ \varphi(e_{j_0}) \neq \sum_{\substack{j \neq j_0 \text, \\ \beta_j \leq \beta_{j_0}}} \lambda_j X^{\beta_{j_0} - \beta_j} \varphi(e_j) \text. \]
    \end{enumerate}

    \begin{proof}
        \ref{prop:reduce:1}: It suffices to check that $\Im \tilde\varphi_{j_0} = \Im \varphi$.
        It is clear that $\Im \tilde\varphi_{j_0} \subseteq \Im \varphi$ and that $\varphi(e_j) \in \Im \tilde\varphi_{j_0}$ for all $j \neq j_0$.
        From \eqref{eqn:minimality3} it follows that $\varphi(e_{j_0}) \in \Im \tilde\varphi_{j_0}$.
        Thus, $\Im \varphi \subseteq \Im \tilde\varphi_{j_0}$ since $\varphi(e_j) \in \Im \varphi_{j_0}$ for all $j \in \{ 1, \dotsc, s \}$.

        \ref{prop:reduce:2}:
        Suppose that there exist $\lambda_1, \dotsc, \lambda_{j_0-1}, \lambda_{j_0+1}, \dotsc, \lambda_s \in k$ such that
        \[ \varphi(e_{j_0}) = \sum_{\substack{j \neq j_0 \text, \\ \beta_j \leq \beta_{j_0}}} \lambda_j X^{\beta_{j_0} - \beta_j} \varphi(e_j) \text. \]
        Then, by Theorem~\ref{thm:minimality2}~\ref{thm:minimality2:6}, the free-injective presentation $\varphi$ is not generator-minimal.

        Conversely, suppose that $\varphi\colon F \to E$ is not generator-minimal.
        Then, according to Corollary~\ref{cor:minimality}~\ref{cor:minimality:3}, there exist homogeneous $\sigma_1, \dotsc, \sigma_s \in R$ not all contained in $\mathfrak m$ such that
        \[ \varphi(\sigma_1 e_1 + \dotsb + \sigma_s e_s) = 0 \text. \]

        Say that $\sigma_{j_0} \in R \setminus \mathfrak m$, i.e.\ $\sigma_{j_0}$ is a monomial of degree zero.
        Then, by subtracting $\varphi(e_{j_0})$ and multiplying with $\sigma^{-1}_{j_0}$ we obtain
        \[ \varphi(e_{j_0}) = - \frac1{\sigma_{j_0}}\sum_{j \neq j_0} \sigma_j \varphi(e_j) \text. \]
        
        Thus, since $\sigma_j \varphi(e_j)$ and $\varphi(e_{j_0})$ are of degree $\beta_{j_0}$, it must hold $\sigma_j = \lambda_j X^{\beta_{j_0} - \beta_j}$ for some $\lambda_j \in k$ if $\beta_{j_0} \geq \beta_j$, and otherwise it must hold $\sigma_j = 0$.
        Multipliying $\lambda_j$ with $-\sigma^{-1}_{j_0}$ and subtituting this in the above equality gives the desired existence result.
    \end{proof}
\end{proposition}

\begin{example}
    We illustrate the application of Proposition~\ref{prop:reduce} to an indecomposable $k[X_1, X_2]$-module, which we take from \cite[Section~3.3]{Lenzen2024}.
    Consider the $k[X_1, X_2]$-graded module~$M$ concentrated in the interval $[(0, 0), (2, 1)]$ that is given by the diagram
    \begin{equation}%
        \label{eqn:example}
        \begin{tikzcd}[ampersand replacement=\&]
            k \& k^2 \& k \\
            0 \& k \& k \text.
            \arrow[from=1-1, to=1-2, "{\left(\begin{smallmatrix}0 \\ 1 \end{smallmatrix}\right)}"]
            \arrow[from=1-2, to=1-3, "{\left(\begin{smallmatrix}1 & 1 \end{smallmatrix}\right)}"]
            \arrow[from=2-1, to=1-1]
            \arrow[from=2-1, to=2-2]
            \arrow[from=2-2, to=1-2, "{\left(\begin{smallmatrix}1 \\ 0\end{smallmatrix}\right)}"]
            \arrow[from=2-2, to=2-3, "1"]
            \arrow[from=2-3, to=1-3, "1"]
        \end{tikzcd}
    \end{equation}

    The monomial matrix of the associated free-injective presentation~$\varphi_M\colon F_M \to E_M$ of~$M$ is given by
    \[ A = \begin{blockarray}{ccccccc}
        & X_1 & X_2 & X_1^2 & X_1X_2 & X_1X_2 & X_1^2X_2 \\
        \begin{block}{c[cccccc]}
            X_1 & 1 & 0 & 0 & 0 & 0 & 0 \\
            X_2 & 0 & 1 & 0 & 0 & 0 & 0 \\
            X_1^2 & 1 & 0 & 1 & 0 & 0 & 0 \\
            X_1X_2 & 1 & 0 & 0 & 1 & 0 & 0 \\
            X_1X_2 & 0 & 1 & 0 & 0 & 1 & 0 \\
            X_1^2X_2 & 1 & 1 & 1 & 1 & 1 & 1 \\
        \end{block}
    \end{blockarray} \text, \]
    where we denote the generator and cogenerator degrees as Laurent monomials, i.e.\ $X^{\alpha}$ denotes the degree $\alpha \in \mathbb Z^n$.
    Let $\mathcal B' = (e_{(1, 0),1}, e_{(0, 1),1}, e_{(2, 0),1}, e_{(1, 1),1}, e_{(1, 1),2}, e_{(2, 1)})$ be the canonical basis of $F_M$.
    Then, by solving the $k$-linear systems induced by Proposition~\ref{prop:reduce}~\ref{prop:reduce:1}, we successively find the relations
    \begin{align*}
        \varphi_M(e_{(2, 0),1}) &= X_1 \cdot \varphi_M(e_{(1, 0),1}) \text, \\
        \varphi_M(e_{(1, 1),1}) &= X_2 \cdot \varphi_M(e_{(1, 0),1}) \text, \\
        \varphi_M(e_{(1, 1),2}) &= X_1 \cdot \varphi_M(e_{(0, 1),1}) \text, \\
        \varphi_M(e_{(2, 0),1}) &= X_1 X_2 \cdot \varphi_M(e_{(1, 0),1}) \text.
    \end{align*}
    Thus, a generator-minimal free-injective presentation of $M$ is given, according to Proposition~\ref{prop:reduce}~\ref{prop:reduce:2}, by the monomial matrix
    \[ \hat A = \begin{blockarray}{ccc}
        & X_1 & X_2 \\
        \begin{block}{c[cc]}
            X_1 & 1 & 0 \\
            X_2 & 0 & 1 \\
            X_1^2 & 1 & 0 \\
            X_1X_2 & 1 & 0 \\
            X_1X_2 & 0 & 1 \\
            X_1^2X_2 & 1 & 1 \\
        \end{block}
    \end{blockarray} \text. \]

    By switching to the Matlis-dual free-injective presentation
    \[ (\hat A)^\vee = \begin{blockarray}{ccccccc}
        & X_1^{-1} & X_2^{-1} & X_1^{-2} & X_1^{-1}X_2^{-1} & X_1^{-1}X_2^{-1} & X_1^{-2}X_2^{-1} \\
        \begin{block}{c[cccccc]}
            X_1^{-1} & 1 & 0 & 1 & 1 & 0 & 1 \\
            X_2^{-1} & 0 & 1 & 0 & 0 & 1 & 1 \\
        \end{block}
    \end{blockarray} \text, \]
    and applying the same procedure to remove redundant cogenerator degrees of $\varphi_M$, we obtain a minimal free-injective presentation $\tilde\varphi$ for~$M$, which is represented by the monomial matrix
    \[ \tilde A = \begin{blockarray}{ccc}
        & X_1 & X_2 \\
        \begin{block}{c[cc]}
            X_1X_2 & 1 & 0 \\
            X_1^2X_2 & 1 & 1 \\
        \end{block}
    \end{blockarray} \text. \]
\end{example}

With Proposition~\ref{prop:reduce} and the previous example in mind, it is possible to devise a technique to iteratively reduce the generators of a given free-injective presentation.
We shall describe the resulting algorithm in the following.
Before that, we explain some notation that is required for the algorithm.

Let $r, s \in \mathbb N$.
For $A \in k^{r \times s}$ and $b \in k^{r}$, we define
\[ L(A, b) = \{ x \in k^s \mid Ax = b \} \text, \]
the \emph{set of solutions of $Ax = b$}.
Note that $L(A, b)$ is either empty or an affine subspace of $k^s$.

Let $A$ be a monomial matrix with cogenerator degrees $(\alpha_i)_{i \in I}$ and generator degrees $(\beta_j)_{j \in J}$.
The \emph{truncation of $A$ at $j_0 \in J$} is the matrix $\trunc(A, j_0) = (c_{ij'})_{i, j' \in J'}$ with $J' = \{ j \in J \mid \beta_j \leq \beta_{j_0} \} \setminus \{ j_0 \}$ such that
\[ c_{ij'} = \begin{cases*}
    a_{ij'} & if $\beta_{j'} \leq \alpha_i$ \\
    0 & otherwise,
\end{cases*} \qquad \text{for $i \in I$, $j' \in J'$.} \]

That is, the matrix $\trunc(A, j_0)$ is constructed as follows:
We regard each columns of $A$ with generator degree less than or equal $\beta_{j_0}$, except for the $j_0$-th column.
Each such column is regarded as a homogeneous element of the injective $k[X_1, \dotsc, X_n]$-module $E = \prod_{i \in I} \Sigma^{\alpha_i} R^\vee$.

For appropriate $j \neq j_0$, the column $A_j$ is multiplied with $X^{\beta_{j_0} - \beta_j}$ to obtain a homogeneous element of $E$ with degree $\beta_{j_0}$, see Example~\ref{ex:action} for the underlying scalar multiplication.
The resulting columns after multiplication are composed to the $k$-valued matrix $\trunc(A, j_0)$.

Lastly, given a monomial matrix
\[ A = \begin{blockarray}{ccccc}
    \beta_1 & \cdots & \beta_j & \cdots & \beta_s \\
    \begin{block}{[ccccc]}
        A_1 & \cdots & A_j & \cdots & A_s \\
    \end{block}
\end{blockarray} \text, \]
with columns $A_1, \dotsc, A_s$, we write
\[ \begin{blockarray}{ccccc}
    \beta_1 & \cdots & \beta_j & \cdots & \beta_s \\
    \begin{block}{[ccccc]}
        A_1 & \cdots & \widehat{A_j} & \cdots & A_s \\
    \end{block}
\end{blockarray} \]
for the monomial matrix obtained by removing the $j$-th column of $A$.

\begin{alg}[Reduction of generators]\label{alg:reduce1}\strut%
    \begin{tabbing}
        \hspace{25mm} \= \qquad \= \qquad \= \kill
        \> \textsc{Name} \' reduce generators \\
        \> \textsc{Input} \' Monomial matrix $A$ with generator degrees $(\beta_j)_{j \in J}$ \\
        \> \textsc{Output} \' Monomial matrix $\tilde A$ such that $\tilde A$ is generator-minimal \\
        \> \textsc{Procedure} \' For each $j \in J$, do \\
        \> \> If $L(\trunc(A, j), A_j) \neq \emptyset$, then \\
        \> \> \> Assign $A \leftarrow \begin{bmatrix} A_1 & \cdots & \widehat{A_j} & \cdots & A_s \end{bmatrix}$. \\
        \> Return $\tilde A \leftarrow A$.
    \end{tabbing}
\end{alg}

The correctness of Algorithm~\ref{alg:reduce1} follows directly from Proposition~\ref{prop:reduce}:
In every iteration, the monomial matrix $A$ is replaced by a monomial matrix~$A$ with the same image by Proposition id., part~\ref{prop:reduce:1}; the assumption is satisfied whenever $L(\trunc(A, j), A_j) \neq \emptyset$.
In the end, Proposition id., part~\ref{prop:reduce:2}, guarantees that after iterating over every $j \in J$, the resulting monomial matrix corresponds to a minimal free-injective presentation.

Since we are dealing only with reflexive modules in this section, and the ring $R = k[X_1, \dotsc, X_n]$ is Noetherian, we can dualize Algorithm~\ref{alg:reduce1} to reduce cogenerators.
Then, we obtain an algorithm to reduce both generators and cogenerators of a monomial matrix~$A$ as follows.

\begin{alg}[Reduction of generators and cogenerators]\label{alg:reduce2}\strut%
    \begin{tabbing}
        \hspace{25mm} \= \qquad \= \qquad \= \kill
        \> \textsc{Name} \' reduce \\
        \> \textsc{Input} \' Monomial matrix $A$ \\
        \> \textsc{Output} \' Monomial matrix $\tilde A$ such that $\tilde A$ is minimal \\
        \> \textsc{Procedure} \' Assign $A \leftarrow \mathop{\text{reduce generators}}(A)$. \\
        \> Assign $A \leftarrow \mathop{\text{reduce generators}}(A^T)^T$. \\
        \> Return $\tilde A \leftarrow A$.
    \end{tabbing}
\end{alg}

\begin{remark}
    Let $\varphi\colon F \to E$ be a free-injective presentation.
    The described procedure effectively determines a set of generators of the graded module $\Ker \varphi/\mathfrak mF$:
    If $\Ker\varphi/\mathfrak mF \cong 0$, then $\Ker\varphi \subseteq \mathfrak mF$, which means that $\varphi$ is generator-minimal.
    
    The Algorithm~\ref{alg:reduce1} can be extended to determine a finite set of generators of the graded module $\Ker  \varphi$, thereby determining the underlying syzygies of a given module\footnote{If $F\overset{q}\twoheadrightarrow M \overset{j}\hookrightarrow E$ is an epi-mono factorization of $\varphi\colon F\to E$ of $M$, then $\Ker\varphi=\Ker q$, i.e.~$\Ker\varphi$ coincides with the syzygy module $\Ker q$ of $M$ with respect to $q\colon F\twoheadrightarrow M$.}.
\end{remark}

\section*{Acknowledgements}

The authors express their gratitude to Ezra Miller for his detailed feedback and comments on an earlier version of this work. 
% \textcolor{blue}{(mention this to the other paper cofiltrarion)} 
% AS gratefully thanks D.~Kozlov for helpful discussions.
The authors would also like to thank Michael Lesnick and Fabian Lenzen for their comments and discussions.

\printbibliography
\end{document}